\documentclass[12pt]{article}
\usepackage[dvips]{graphicx}
\usepackage{latexsym}
\usepackage{amssymb}

\newcommand{\C}{\mathbb{C}}
\newcommand{\CP}{\mathbb{CP}}
\newcommand{\RP}{\mathbb{RP}}
\newcommand{\spp}{\mathbb{S}}
\newcommand{\R}{\mathbb{R}}

\newcommand{\PP}{\mathbb{P}}

\renewcommand{\d}{\mathrm{d}}

\newcommand{\koniec}{\begin{flushright}  $\Box $ \end{flushright}}
\def\be{\begin{equation}}
\def\ee{\end{equation}}
\def\theequation{\thesection.\arabic{equation}}

\def\Sm{\Sigma}

\def\O{\cal O}
\def\om{\omega}

\def\p{\partial}

\def\ll{\lambda}

\def\a{\alpha}
\def\l{\lambda}

\def\O{{\cal O}}

\def\z{\cal Z}
\newtheorem{theo}{Theorem}[section] 
  
\newtheorem{lemma}[theo]{Lemma}
\newtheorem{defi}[theo]{Definition}

\begin{document}
\pagestyle{plain}

\title{\vskip -70pt
\begin{flushright}
{\normalsize DAMTP-2005-20} \\
\end{flushright}
\vskip 80pt
{\bf Paraconformal geometry of  $n$th order ODEs, and exotic
holonomy in dimension four}\vskip 20pt}

\author{Maciej Dunajski\thanks{email M.Dunajski@damtp.cam.ac.uk}\\[15pt]
{\sl Department of Applied Mathematics and Theoretical Physics} \\[5pt]
{\sl University of Cambridge} \\[5pt]
{\sl Wilberforce Road, Cambridge CB3 0WA, UK} \\[15pt]
and \\[10pt]
Paul Tod\thanks{email tod@maths.ox.ac.uk} \\[15pt]
{\sl The Mathematical Institute}\\[5pt]
{\sl Oxford University}\\[5pt]
{\sl 24-29 St Giles, Oxford OX1 3LB, UK}\\[15pt]
}
\date{} 
\maketitle
\begin{abstract}
We characterise $n$th order ODEs for which the space of solutions $M$
is equipped with a particular 
paraconformal structure in the sense of \cite{BE}, 
that is a splitting of
the tangent bundle as a symmetric tensor product of rank-two vector bundles. 
This leads to the 
vanishing of $(n-2)$ quantities constructed from of the ODE.

If $n=4$ the paraconformal structure is shown to be equivalent to
the exotic ${\cal G}_3$ holonomy of Bryant.
If $n=4$, or $n\geq 6$ and  $M$ admits 
a  torsion--free connection compatible with the
paraconformal structure then the ODE is trivialisable by point or contact 
transformations respectively. 

If $n=2$ or $3$  $M$ admits an  affine paraconformal 
connection with no torsion. In these cases 
additional constraints can be imposed on
the ODE so that $M$ admits a projective structure if $n=2$, 
or an Einstein--Weyl structure if $n=3$. 
The third order ODE can in this case be 
reconstructed from the Einstein--Weyl data.
\end{abstract}
\newpage
\section{Introduction}
\setcounter{equation}{0}
Consider a relation of the form
\be
\label{relation}
\Psi(x, y, {t})=0
\ee
between the variables ${t}= (t^1, t^2, ..., t^n)$ (local coordinates
on an $n$-dimensional manifold $M$), and $(x, y)$ (local coordinates on a 
two-dimensional  manifold ${\z}$, which we shall call 
the twistor space). 
For each fixed choice of $(x, y)$m the relation (\ref{relation}) 
defines a hypersurface in  $M$. Conversely each choice of ${t}$ defines 
a curve $L_t$ in ${\cal Z}$. Given conditions on the derivatives of $\Psi$, 
we can  apply the implicit function theorem to 
(\ref{relation}), and regard $L_t$ as a graph
\be
\label{graph}
x\longrightarrow (x, y=Z(x, {t})).
\ee
Consider the system of algebraic equations consisting of
$y=Z(x, {t})$, and the first $(n-1)$ derivatives with respect to $x$.
Solving this system for ${t}$, and differentiating once more with respect
to $x$ yields
\be
\label{ode}
y^{(n)}:=\frac{\d^n y}{\d x^n}=F(x, y, y', ..., y^{(n-1)}),
\ee 
where the explicit form of $F$ is completely determined by (\ref{relation}).
This procedure will lead to an  $n$th (as opposed to a lower order) order
ODE if $\Psi$ is sufficiently smooth and  non-degenerate in an suitable sense. 
This open  non-degeneracy condition is
best expressed in terms of $Z(x, t)$ by demanding that the gradients
$\nabla Z, \nabla Z', ..., \nabla Z^{(n-1)}$ with respect to $t$ 
are linearly independent
on $M$.

To achieve a more geometric picture we define
the $(n+1)$-dimensional hyper-surface ${\cal F}\subset {\cal Z}\times M$ 
called the correspondence space by 
the following incidence relation
\be
\label{correspondence}
{\cal F}=\{ ((x, y), {t}) \in {\cal Z}\times M| z\in L_t\}.
\ee
The double fibration
\be
\label{doublefib}
{M}\stackrel{p}\longleftarrow 
{\cal F}\stackrel{q}\longrightarrow {\z}
\ee
is then defined by the relation (\ref{relation}), and therefore by the ODE 
(\ref{ode}).

Putting various geometric structures on $M$ (which from now on will be 
identified with the space of solutions to the ODE (\ref{ode})) imposes 
additional constraints on  $F$. This idea goes back to Cartan
\cite{C41}, 
and his
program of `geometrising' ODEs. 
Extending Cartan's program to PDEs is possible, and underlies
some approaches to general relativity \cite{FKN01}, and 
other problems in mathematical physics \cite{Druma}.

A different approach based on twistor theory was suggested by Hitchin 
\cite{H82}. In this approach one works in the holomorphic category and 
$(x, y, {t})$ are complex numbers.
The graph (\ref{graph}) represents
 a compact holomorphic (i.e. rational)
curve in $\z$ with a prescribed normal bundle. The local differential 
geometry of $M$ is encoded in  global embeddings of the family of curves 
(parametrized by $t$) in $\z$.
In this approach, 
the ODE (\ref{ode}) does not explicitly appear in the correspondence
between $M$ and $\z$. The details of Hitchin's construction and its connection 
with the ODE approach have partially been worked out only for $n=2$ 
\cite{Mc05}. 
In this case there exists an embedding of a rational curve
with  normal bundle $\O(1)$ in $\z$ if and only if 
\be
\label{geodesiccond}
\frac{\d^2}{\d x^2}F_{11}-4\frac{\d}{\d x}F_{01}-F_1\frac{\d}{\d x}F_{11}
+4F_1F_{01}-3F_0F_{11}+6F_{00}=0,
\ee
where 
\[
F_0=\frac{\p F}{\p y}, \qquad F_1=\frac{\p F}{\p y'}, \qquad F_2=\frac{\p F}{\p y''},\qquad 
..., \qquad
F_{n-1}=\frac{\p F}{\p y^{(n-1)}},
\]
and
\[
\frac{\d}{\d x}=\frac{\p }{\p x}+
\sum_{k=1}^{n-1}y^{(k)}\frac{\p }{\p y^{(k-1)}}+F\frac{\p}{\p y^{(n-1)}}.
\]
The
two-dimensional moduli
space $M$ of $\O(1)$ curves  is 
in this case equipped with a projective structure, in the sense that
the hyper-surfaces (curves) of constant $(x, y)$ in (\ref{relation})
are geodesics of a torsion-free connection. Conversely, given a
projective structure on $M$ one defines ${\cal Z}$ as the quotient space
of the foliation of  $\PP (TM)$ by the orbits of the geodesic flow. 
Each projective tangent space $\PP(T_{t}M)$ maps to a rational
curve with self-intersection number one in ${\cal Z}$.

The case $n=3$ goes back to Cartan \cite{C41} and Chern
\cite{Ch40}, and  was recently revisited in \cite{T00}. 
The conformal structure on $M$ is defined by demanding
that hyper-surfaces $\Sm\subset M$ corresponding to points in $\z$ 
are null. This conformal structure is well defined
if $F(x, y, y', y'')$ satisfies a third--order differential
constraint 
\be
\label{wun3}
\frac{1}{3}F_{{2}}{\frac {\d}{\d x}}F_{{2}}-
\frac{1}{6}{\frac {\d^{2}}{\d{x}^{2}}
}F_{{2}}+\frac{1}{2}{\frac {\d}{\d x}}F_{{1}}-{\frac {2}{27}}\left (F_{
{2}}\right )^{3}-\frac{1}{3}F_{{2}}F_{{1}}-F_{{0}}=0.
\ee
This constraint has already appeared in the work of
W\"unschmann \cite{W05}. 
The hyper-surfaces $\Sm$ are totally geodesic with respect to some 
torsion-free connection $D$ if $F$ satisfies the additional condition 
\be
\label{EWcond}
\frac{\d^2}{\d x^2}F_{22}-\frac{\d}{\d x}F_{12}+F_{02}=0.
\ee
The existence of a two-parameter family
of totally geodesic null hypersurface in $M$ is equivalent to the
vanishing of the trace-free part of the symmetrised Ricci tensor of $D$.
This is the Einstein--Weyl condition first introduced in \cite{C43}. 
The three-dimensional Einstein--Weyl spaces can therefore be obtained
from a particular class of third-order ODEs (\ref{ode}).
In the twistor approach \cite{H82} the moduli space of 
rational curves in $\z$ with  normal bundle $\O(2)$ is automatically
equipped with an EW structure, and all analytic 
EW structures locally arise in such a way.

 The only other case which has attracted some attention is $n=4$.
Bryant \cite{B91} has shown that there exists a correspondence between 
a class of fourth order ODEs, and exotic non-metric holonomies in dimension 
four. The conditions on $F$ are only implicit in Bryant's work.

In Section \ref{sec_wun} we shall generalise the W\"unschmann
condition (\ref{wun3}) to $(n-2)$ conditions in the case of $n$th order
ODEs, and give an example of an ODE for which all these conditions
are satisfied. 
\begin{defi}
A paraconformal structure on a smooth manifold $M$ is a bundle isomorphism
\be
\label{paracon}
TM\cong \spp \odot \spp \odot ... \odot \spp =\mbox{S}^{n-1}{(\spp)},
\ee
where $\spp\rightarrow M$ is a real rank--two vector bundle, and $\odot$
denotes symmetric tensor product.
\end{defi}
More general paraconformal structures have been considered 
in \cite{BE} and \cite{AG98} (where they were called almost 
Grassmann structures)
but we shall only work with (\ref{paracon}).
The isomorphism (\ref{paracon}) identifies each tangent space $T_{t} M$  
with the space of homogeneous 
$(n-1)$th order polynomials 
in two variables. The vectors corresponding to polynomials with repeated root of multiplicity $(n-1)$ are called maximally null. A hypersurface in $M$
is maximally null if its normal vector is maximally null. 
All results established in this paper are local on $M$.

In the next section we shall prove the following
\begin{theo} 
\label{theo_wun}
Assume that the space of solutions $M$ to the $n$th order ODE (\ref{ode}) is
equipped with a paraconformal structure (\ref{paracon}) such that the 
two--parameter family of hypersurfaces (\ref{relation}) are maximally null.
Then $F$ satisfies $(n-2)$ 
conditions of the form
\be
\label{invariantsC}
C_k\Big(F_i, \frac{\d F_i}{\d x}, ..., 
\frac{\d^{n-1} F_i}{\d x^{n-1}}\Big)=0, \qquad i=0, ..., n-1, 
\qquad k=1, ..., n-2.
\ee
Each expression $C_k$ is a polynomial in the derivatives of $F$ of order
less than or equal to $n$, and equations (\ref{invariantsC})
are invariant under point transformations on ${\cal Z}$ (i.e. transformations
induced by a change of variables $\hat{x}=\hat{x}(x, y),
\hat{y}=\hat{y}(x, y)$).

 Conversely given (\ref{invariantsC}) there is a paraconformal structure on $M$ such that points in $\z$ define maximally null hypersurfaces in $M$.
\end{theo}

In Section \ref{secdoubrov} we shall relate the paraconformal structure
on $M$ to the existence of a subbundle of $\PP(TM)$ with rational normal
curves as fibres, and demonstrate that the conditions (\ref{invariantsC})
are the classical Wilczynski invariants \cite{Wilczynski} on the linearisation
of (\ref{ode}) at a fixed solution. This strengthens Theorem 
(\ref{theo_wun}) as the Wilczynski conditions are invariant under
the wider class of contact transformation. 

In Section \ref{sec_twistor} we shall show that
the paraconformal structure (\ref{paracon}) exists on the moduli
space of rational curves with normal bundle ${\cal O}(n-1)$ in a
complex surface, which leads to a twistorial interpretation of the 
constraints on $F$. 

In Section \ref{sec_holo} we shall discuss the case   
$n=4$, where  the existence of the paraconformal structure is 
equivalent to the existence of the torsion--free connection 
with  ${\cal G}_3$ holonomy on the space of solutions to
(\ref{ode}).
\begin{theo}
\label{theo_G3}
Let $M$ be the space of solutions to the fourth order
ODE
\[
\frac{\d^4 y}{\d x^4}=F(x, y, y', y'', y''').
\]
The following conditions are equivalent:
\begin{enumerate} 
\item $M$ admits the paraconformal structure (\ref{paracon}) with maximally 
null surfaces (\ref{relation}). 
\item $M$ admits a torsion--free connection with holonomy ${\cal G}_3$
\item $F$ satisfies a pair of third order PDEs
\be
\label{W4_1}
{\frac {11}{1600}}\,\left (F_{{3}}\right )^{4}-{\frac {9}{50}}\,
\left (F_{{3}}\right )^{2}{\frac {\d}{\d x}}F_{{3}}-{\frac {1}{200}
}\,\left (F_{{3}}\right )^{2}F_{{2}}+{\frac {21}{100}}\,\left ({
\frac {\d}{\d x}}F_{{3}}\right )^{2}+{\frac {1}{50}}\left ({\frac {\d
}{\d x}}F_{{3}}\right )F_{{2}}
\ee
\[
-{\frac {9}{100}}\,\left (F_{{2}}
\right )^{2}+{\frac {7}{20}}\,F_{{3}}{\frac {\d^{2}}{\d{x}^{2}}}F_{{3
}}-\frac{1}{5}\,{\frac {\d^{3}}{\d{x}^{3}}}F_{{3}}+\frac{3}{10}\,{\frac {\d^{2}}{\d{
x}^{2}}}F_{{2}}-\frac{1}{4}\,F_{{3}}{\frac {\d}{\d x}}F_{{2}}-F_{{0}}=0,
\]
\be
\label{W4_2}
\frac{9}{4}\,F_{{3}}{\frac {\d}{\d x}}F_{{3}}-\frac{3}{2}\,{\frac {\d^{2}}{\d{x}^{2}}
}F_{{3}}+3\,{\frac {\d}{\d x}}F_{{2}}-\frac{3}{8}\,\left (F_{{3}}\right 
)^{3}-\frac{3}{2}\,F_{{2}}F_{{3}}-3\,F_{{1}}=0.
\ee
\end{enumerate}
(differentiating the second condition w.r.t $x$ and subtracting its
constant multiple from the  first one leads to a couple of $3rd$ order
PDEs for $F$). 
\end{theo}
In Section \ref{sec_conn} we shall study connections preserving the
paraconformal structure, and show that they must necessarily have
torsion if $n=4$ or $n\geq 6$ 
(we stress that our definition of the torsion--free
paraconformal connection is stronger that Bryant's torsion--free
${\cal G}_3$ holonomy \cite{B91}).
\begin{theo}
\label{connection}
If  $M$ admits the paraconformal structure (\ref{paracon}) and 
\[
D:\Gamma(\spp^k)\longrightarrow \Gamma(\spp^k\otimes T^*M)=
\Gamma(\spp^{k+n-1}),
\] 
where $\spp^k=\spp^{\otimes k},\;\;k=0, 1, ..., n-1$,
is a torsion--free connection preserving the paraconformal structure
then additional constraints {\em(\ref{other_cond}, \ref{A1}--\ref{A5})},
with $(P, Q)$ given by {\em(\ref{pandq})},  
need to be satisfied.

In particular
the ODE {\em(\ref{ode})} can be transformed to 
\[
\frac{\d ^n y}{\d x^n}=0
\]
by point transformations if
$n=4$, and by contact transformations if $n\geq 6$.
If $n<4$ then $D$ always exists.
\end{theo}
The case $n=5$ is special, and we shall give an example 
(equation (\ref{order_5}))
of a 5th order ODE such that 
$M$ admits a paraconformal structure with a paraconformal 
torsion--free connection.

If $n=3$
the paraconformal structure (\ref{paracon}) is conformal in the usual
sense. The rank--two vector bundle $\spp$ in the isomorphism 
(\ref{paracon}) is in this case the usual spin bundle and its sections are 
called the two--component spinors. We shall continue to call it the spin bundle
in the general case $n>3$ 
although its sections are not spinors as there is no
underlying orthogonal group.

In Section \ref{sec_EW} we shall concentrate on the case
$n=3$, where  the condition  (\ref{invariantsC}) 
for $F$ is the W\"unschmann condition
(\ref{wun3}). We shall give an algorithm 
for determining the third order ODE satisfying (\ref{EWcond}), and
(\ref{wun3}) from a given
Einstein--Weyl structure, based on the Lax formulation of the
Einstein--Weyl conditions. 

Most calculations leading  to invariants like (\ref{W4_1}, \ref{W4_2})
were performed (or checked) using MAPLE. The resulting long expressions are
usually unilluminating. They are nevertheless useful in constructing
explicit examples like (\ref{example1}), and we have decided to include them in
the paper. Readers who want to verify our calculations can obtain
the MAPLE programs from us.
\section{W\"unschmann invariants}
\setcounter{equation}{0}
\label{sec_wun}
In this Section we shall establish Theorem \ref{theo_wun}, and give an
example of an ODE which leads to a paraconformal structure for any $n$.
First we need to introduce some notation.

The isomorphism (\ref{paracon}) identifies each tangent space $T_{t} M$  with the space of homogeneous 
$(n-1)$th order polynomials 
in two variables
\[
T\in TM\longrightarrow {\bf t}=T^{A_1A_2...A_{n-1}}z_{A_1}z_{A_2}...z_{A_{n-1}},\qquad
A_1, A_2, ..., A_{n-1}=0, 1
\]
where $z_{A_i}=(z_0, z_1)\in\R^2$, and $T^{A_1A_2...A_{n-1}}$ is symmetric
in its indices.

Our considerations are local on $M$ so we choose a trivialisation of $TM$
and represent a vector $T$ by its components $T^a, a=1, ..., n$ with respect  
to some basis. We also choose a trivalisation of $\spp$. The paraconformal
structure is then defined in terms of van der Waerden symbols 
$\sigma^a_{A_1...A_{n-1}}$ (which are symmetric in $A_1, ..., A_{n-1}$) by
\[
T^a=\sigma^a_{A_1...A_{n-1}} T^{A_1A_2...A_{n-1}}.
\] 
The bold letters denote  homogeneous polynomials.
The summation convention is used unless stated
otherwise.
For any $0\leq r\leq n-1$
we define  $V_{r}\subset \R[z_0, z_1]$ to be the $(r+1)$-dimensional space of homogeneous polynomials 
of degree $r$. Let ${\bf t}\in V_{n-1}$. 
The space $V_{n-1}$ is an $SL(2, \R)$ module, and the infinitesimal action 
of $SL(2, \R)$ is generated by ${\bf t}\rightarrow H({\bf t})$, where
\[
H={H_A}^Bz_B\frac{\p}{\p z_A} \in {\bf sl}(2, \R),
\]
and ${H_A}^B$ is one of the following matrices
\[
\left (
\begin{array}{cc}
0&1\\
0&0
\end{array}
\right ),\qquad
\left (
\begin{array}{cc}
0&0\\
1&0
\end{array}
\right ),\qquad
\left (
\begin{array}{cc}
1&0\\
0&-1
\end{array}
\right ).
\]
For each $p\geq 0$ define a linear, $SL(2, \R)$ equivariant 
mapping $V_r\otimes V_s\longrightarrow V_{r+s-2p}$ given by
\be
\label{inner_product}
<{\bf t}, {\bf s}>_p
=\varepsilon_{A_1B_1}\varepsilon_{A_2B_2}...
\varepsilon_{A_pB_p}\frac{\p^p {\bf t}}{\p z_{A_1}...\p z_{A_p}}
\frac{\p^p {\bf s}}{\p z_{B_1}...\p z_{B_p}},
\ee
where
\[
\varepsilon_{AB}=\left (
\begin{array}{cc}
0&1\\
-1&0
\end{array}
\right )
\]
is a symplectic form on the real fibres of $\spp$. In the following sections 
we will not fix this  symplectic form.
The anti--symmetric matrix
$\varepsilon$ will only be defined up to scale, and each choice of
the scale will provide an identification between $\spp$ and its dual bundle.

In particular $< , >_{n-1}:V_{n-1}\times V_{n-1}\rightarrow \R$ is a symmetric
or skew-symmetric (depending on $n$) bilinear form on $V_{n-1}$.
For $m=0, 1, ..., n-1$ define $C_m\subset V_{n-1}$ to be a two-dimensional {\em cone of
order m}, given by all polynomials ${\bf t}={\bf p}^m {\bf r}$,
where ${\bf p}\in V_1$, and ${\bf r}\in V_{n-m-1}$.

{\bf Proof of Theorem \ref{theo_wun}}.
We want  to define a paraconformal structure by requiring 
\[
\sigma^a_{A_1...A_{n-1}}\frac{\p Z}{\p t^a}=p_{A_1}...p_{A_n},
\]
where $y=Z(x, t)$ is a surface in $M$  corresponding to a solution
of $(\ref{ode})$. The symbols $\sigma$ are of course independent of $x$,
but the $p_A$s (and the corresponding polynomials) do depend on $x$ as does
$Z$. 

We shall first assume that this paraconformal structure exists
on $M$ and establish the conditions satisfied by the ODE.

All sections  of $\spp$ correspond to degree one
homogeneous polynomials in $\zeta^A$. 
If $(p_0, p_1)\in \R^2$, then ${\bf p}=p_0\zeta^0+p_1\zeta^1\in V_1$. 
Let ${\bf p}\in V_1$, and let
${\bf q}={\bf p}'$. The fibres of the spin bundle are two-dimensional therefore
\be
\label{2ds}
\frac{\d \bf q}{\d x}=P{\bf p}+Q{\bf q}
\ee
for some $P, Q\in V_0$. These relations hold on the correspondence space
(\ref{correspondence}) and ${\bf p}, {\bf q}$ are regarded as $x$--dependent 
sections of the 
bundle $\spp$  pulled back from $M$ to ${\cal F}$. 
 
Consider an element $T$ of $C_{n-1}$,  (a maximally null vector). 
The maximally null vectors correspond to polynomials ${\bf T}={\bf
p}^{n-1}$ with a repeated root of
multiplicity $(n-1)$.
Each solution $y=Z(x, {t})$ defines a section of $T^*M$ given by the gradient
on $M$
\[
E=\nabla Z.
\]
Assume that $E$ is a maximally null one-form (an element of $C_{n-1}$), 
and  construct a frame of $n$ one-forms on $T^*M$ given by 
\[
{ E}, {E}', { E}'',..., { E}^{(n-1)}. 
\]
Requiring that 
these forms are linearly independent imposes a condition of 
non-degeneracy, which is an open condition
on $\Psi$ in the relation (\ref{relation}). 

The corresponding polynomials are of the form
\begin{eqnarray*}
{\bf E}&=& {\bf p}^{n-1}\\
{\bf E}'&=&0+a_{11}{\bf p}^{n-2}{\bf q}\\
{\bf E}''&=&a_{20} {\bf p}^{n-1}+a_{21}{\bf p}^{n-2}{\bf q}+
a_{22}{\bf p}^{(n-3)}{\bf q}^2\\
&.&\\
&.&\\
&.&\\
{\bf E}^{(n-1)}&=&a_{(n-1)0} {\bf p}^{n-1}+a_{(n-1)1}{\bf p}^{n-2}{\bf q}
+...+a_{(n-1)(n-1)}{\bf q}^{n-1},
\end{eqnarray*}
or in general
\be
\label{in_general}
{\bf E}^{(i)}=\sum_{k=0}^i a_{ik}{\bf p}^{n-1-k}{\bf q}^k, \qquad i=0,..., n-1.
\ee
The upper triangular matrix $(a_{ij})$ can be computed using
(\ref{2ds}). It  depends on $P, Q$ and their derivatives with
respect to $x$. 
The $n$th derivative of ${\bf E}$ with respect to $x$ is given by
\[
{\bf E}^{(n)}=a_{n0} {\bf p}^{n-1}+a_{n1}{\bf p}^{n-2}{\bf q}
+...+ a_{n(n-1)}{\bf q}^{n-1}.
\]
Remembering that  $y=Z(x, {t})$ is a solution to (\ref{ode}), and 
using the chain rule we  express ${\bf E}^{(n)}$
as a linear combination of ${\bf E}, {\bf E}', ..., {\bf E}^{(n-1)}$ by
\be
\label{expE}
{\bf E}^{(n)}=\sum_{i=0}^{n-1}F_i {\bf E}^{(i)}.
\ee
This gives rise to 
\be
\label{matrixODE}
a_{nj}=\sum_{i=0}^{n-1}F_i a_{ij}.
\ee
Solving  these $n$ equations for $P$ and $Q$ yields
\begin{eqnarray}
\label{pandq}
Q&=&\frac{2}{n(n-1)}F_{n-1}\nonumber\\
P&=&\frac{1}{n(n^2-1)}\Big(\frac{(3n-1)(n-2)}{n(n-1)}F_{n-1}^2+6F_{n-2}
-2(n-2)\frac{\d}{\d x}F_{n-1}\Big)
\end{eqnarray}
(the calculations leading to these formulae are presented in the Appendix).
The remaining equations imply the vanishing of
$(n-2)$ expressions constructed out of  
$F$. Each expression is a polynomial $C_k$ in the derivatives of $F$ of order
less than or equal to $n$, of the form (\ref{invariantsC}).
The vanishing of these expressions  characterises a class of ODEs (\ref{ode}) such that
their solution spaces admit a paraconformal structure (\ref{paracon}).
It follows from the construction (and it
may be checked if desired) that vanishing of these expressions is
invariant under point transformations
\[
\hat{x}=\hat{x}(x, y ),\qquad \hat{y}=\hat{y}(x, y).
\]

Conversely, let us assume that we are given an ODE (\ref{ode}) such that
the conditions (\ref{invariantsC}) hold. We define $P$ and $Q$ by
(\ref{pandq}) and solve the linear system of ODEs 
${\bf q}={\bf p}', {\bf q'}=P{\bf p}+Q{\bf q}$ to determine 
${\bf p}$ and ${\bf q}$. We then define a basis of $TM$ by
the procedure of taking gradients leading to  
(\ref{in_general}). The consistency conditions are guaranteed by the 
calculation in the first part of the proof. This gives a paraconformal 
structure such that the surfaces $y=Z(x, t)$ in $M$ are totally null.

\koniec
If $n=2$ the paraconformal structure always exists. If $n=3$ the
condition (\ref{invariantsC}) is given by  (\ref{wun3}).
If $n=4$ we have two conditions given by (\ref{W4_1}) and (\ref{W4_2}).
The polynomials corresponding to $n=5$ are given in the Appendix.

There are no terms of the form $F_{ij}$ 
(recall that $F_{ij}=\p^2 F/\p y^{(i)}\p y^{(j)}$)
in $C_k$, so the coefficients 
of the polynomials $C_k$ are determined by looking at the special case
of linear homogeneous equations, where
\[
F=p_{n-1}(x)y^{(n-1)}+...+p_0(x)y,
\]
and $F_k=p_k(x)$.

\subsection{Example}
The general solution of (\ref{matrixODE}) regarded as an overdetermined
system of PDEs for $F$ appears to be intractable. 
To find some examples we seek
$F=F(y^{(n-1)}, x)$. Let $z:=y^{(n-1)}$.
The ODE (\ref{ode}) reduces to a 1st order ODE, and 
a sequence of quadratures
\[
z'=F(z, x), \qquad y(x)=\int^x\int^{x_{n-1}}...\int^{x_2}z(x_1)\d x_1\d x_2...\d x_{n-1}.
\]
We use MAPLE to verify that 
all constraint equations (\ref{invariantsC}) reduce to
\[
\frac{\d}{\d x}F_{n-1}=\frac{1}{n}(F_{n-1})^2,\qquad\mbox{where now}\;\;
\frac{\d}{\d x}=\frac{\p}{\p x}+F\frac{\p}{\p z}.
\]
We therefore need to solve
\[
\frac{\p^2 F}{\p z\p x}=\frac{1}{n}\Big(\frac{\p F}{\p z}\Big)^2
-F\frac{\p^2 F}{\p z^2}.
\]
The Legendre transformation
$s=F_z, G(s, x)=F(z(s,x),x)-sz(s,x)$ gives a linear equation for $G$.
This transform can then be inverted, and
the solution can be found for an arbitrary `initial data' $F(z, 0)$.
To write down an explicit example make a further assumption that
$F$ is independent of $x$, which yields
\[
F(z)=(az+b)^{\frac{n}{n-1}}.
\]
Redefining $y(x)$ by a point transformation we can set $b=0$,
so the $n$th order ODE is
\be
\label{example1}
y^{(n)}=(a y^{(n-1)} )^{\frac{n}{n-1}}.
\ee
The corresponding $n$-parameter family of solutions 
to (\ref{ode}) is readily found
\[
y=t^1+t^2x+...+t^{n-1}x^{n-2}
-\frac{(n-1)^{(n-1)}}{a^n(n-2)!}\ln{(x+t^{n})}
\]
If $n=3$ the  
space of solutions to (\ref{ode}) is equipped with the 
NIL Einstein--Weyl structure \cite{T00} (see also Section \ref{sec_EW}).
\section{Comparison with Doubrov--Wilczynski invariants}
\setcounter{equation}{0}
\label{secdoubrov}
In this section we shall demonstrate that for a given ODE
the generalised W\"unschmann conditions (\ref{invariantsC}) 
are   equivalent to the vanishing of a set of invariants 
constructed by Doubrov \cite{Dou}. Doubrov's work builds on an old 
theorem of Wilczynski \cite{Wilczynski}, which we shall review first.

Let
\be
\label{linearODE}
y^{(n)}=p_{n-1}(x)y^{(n-1)}+...+p_0(x)y
\ee
be a linear homogeneous nth order ODE, defined up to 
transformations of the form
\[
(y, x)\longrightarrow (a(x)y, b(x)).
\] 
Wilczynski has demonstrated that the ODE (\ref{linearODE}) is 
trivialisable by this transformation 
if $(n-2)$ expressions constructed out of the $p_i(x)$s
and their derivatives vanish. More formally he has shown the following
\begin{theo}[Wilczynski \cite{Wilczynski}]
\label{th_wil}
The following three
conditions are equivalent:
\item (a) 
The Wilczynski invariants 
\be
\label{wilczynski}
\theta_k\Big(p_i, \frac{\d p_i}{\d x}, ..., 
\frac{\d^{n-1} p_i}{\d x^{n-1}}\Big)=0, \qquad i=0, ..., n-1, 
\qquad k=1, ...n-2
\ee
of the  equation {\em(\ref{linearODE})} vanish.
\item (b) The equation {\em(\ref{linearODE})} 
can be transformed to $y^{(n)}=0$ by a 
change of variables 
$(y,x)\rightarrow (a(x)y, b(x))$.
\item (c) Let $y_1(x), ..., y_n(x)$ be any basis of the solution space of
{\em(\ref{linearODE})}. 
Then the embedding $\RP^1\rightarrow \RP^{n-1}$ given by
$x\rightarrow [y_1(x) : ... : y_n(x)] \subset \RP^{n-1}$, is an open subset 
of the normal rational curve.
\end{theo}
The expressions $\theta_k$ are polynomials in $p_i$s and their derivatives of
order less than or equal to $n$. Wilczynski has produced compact 
formulae for these polynomials if $p_{n-1}(x)=p_{n-2}(x)=0$ (it is
always possible to set these two coefficients to zero by a choice of 
$a(x), b(x)$). We shall not need these formulae in what
follows.

The Doubrov invariants introduced in \cite{Dou}  for general 
ODEs are the classical Wilczynski invariants
(\ref{wilczynski}) of 
their linearizations. That is, when restricted to any
solution (\ref{ode}), they produce the classical Wilczynski invariants of the
linearization of (\ref{ode}) at this solution.
In practise, Doubrov's invariants are calculated by replacing
$p_k(x)$ by $F_k(x, y, ..., y^{(n-1)})$ in Wilczynski's invariants.
\begin{theo}[Doubrov \cite{Dou}]
\label{th_doubrov1}
Let $\theta_k$ be the Wilczynski invariants of a linear $n$th order
 homogeneous
ODE. Then the vanishing of all the expressions
\be
\label{doubrov1}
L_k=\theta_k\Big(F_i, \frac{\d F_i}{\d x}, ..., 
\frac{\d^{n-1} F_i}{\d x^{n-1}}\Big),\qquad i=0, ..., n-1, 
\qquad k=1, ...n-2,
\ee
is invariant under
contact transformations 
of the general $n$th order ODE {\em(\ref{ode})}.
\end{theo}
The invariants $L_k$ will be from now on called the Doubrov invariants.

The ideas behind the proof of the next theorem are due to Doubrov. The 
geometric content of this theorem (constraining a curve until
it becomes a rational normal curve) was also  known to Bryant 
\cite{Bprivate}.
\begin{theo} 
\label{wun_doubrov}
The vanishing of the generalised W\"unschmann conditions {\em{(\ref{invariantsC})}}
is equivalent to the vanishing of the Doubrov 
invariants {\em{(\ref{doubrov1})}}.
\end{theo}
To prove this theorem we shall need a Lemma
\begin{lemma}
\label{lemma_rational_normal}
The existence of the paraconformal structure $TM=\spp\odot \spp\odot ... \odot \spp$
is equivalent to the existence of a sub-bundle of $\PP(TM)$ with rational
normal curves as fibres.
\end{lemma}
{\bf Proof.}
For any rational normal curve ${\cal C}\subset \PP(V)$, where V is an n-dimensional vector
space there exists an isomorphism $V\rightarrow V_{n-1}$
which identifies ${\cal C}$  with $x\rightarrow(1, x, ..., x^{n-1})$. 
Any two such isomorphisms differ by a $PGL(2, \R)$ projective
transformation.

Now if $M$ is paraconformal then $T_tM$ is isomorphic to $V_{n-1}$ for all
$t\in M$. Any element ${\bf p}\in \spp$ (a two dimensional spinor) 
gives a totally null cone ${\bf p}^{n-1}\in T_tM$, or a rational 
normal ${\cal C}_t$ in $P(T_tM)$. The union of ${\cal C}_t$ as $t$
varies in $M$ gives a 
sub-bundle of $\PP(TM)$.
\koniec
{\bf Proof of Theorem \ref{wun_doubrov}.}
If the Doubrov invariants (\ref{doubrov1}) 
vanish, then any linearisation of
(\ref{ode}) has vanishing Wilczynski invariants, and by Theorem \ref{th_wil}
we have a
rational normal curve in the projectivization of each tangent space to
the solution space. This, by Lemma \ref{lemma_rational_normal}, 
is equivalent to the existence of 
the paraconformal structure on the solution space. Theorem  \ref{theo_wun}
now implies that  the conditions (\ref{invariantsC}) are satisfied.

Conversely, let the ODE (\ref{ode}) satisfy (\ref{invariantsC}), and let
$y=Z(x, t), t=(t^1, ..., t^n)$ be its general solution. Fix $t=T_0$. 
The vector $E= \p_a y$ in $T_{T_0}M $, where $\p_a=\p/\p t^a$, 
is maximally null with respect to the paraconformal
structure defined by the ODE, so that it is given by ${\bf p}^{n-1}$. 
Each component of $E=\nabla Z$ 
is a derivative of a solution to the ODE (\ref{ode}), 
and as such it satisfies the linearisation of (\ref{ode}) 
around $Y=Z(x, T_0)$, as the differentiation of (\ref{ode}) yields
\be
\label{linear_1}
\frac{\p y^{(n)}}{\p t^a}=
\sum_{i=0}^{n-1} \Big(\frac{\p F}{\p y^{(i)}}|_{y=Y}\Big) 
\frac{\p y^{i}}{\p t^a}.
\ee
These equations are homogeneous linear ODEs for $\p y/\p t^a$.
But then the nullity of $E$ implies
that the embedding $x\rightarrow [\p_1 y (x): ... : \p_n y(x)]$  
is a rational normal curve for any $T_0$, and so 
the linearised ODE (\ref{linear_1}) is trivialisable by 
Wilczynski's theorem \ref{th_wil}. 
Therefore the ODE (\ref{ode}) has vanishing   Doubrov invariants 
(\ref{doubrov1})
which proves the relative equivalence of the two sets of invariants.
(i.e their vanishing is equivalent).
\koniec
In Section \ref{sec_conn} we shall need another result of Doubrov's
\begin{theo}[Doubrov \cite{Dou}]
\label{doubrov2}
The ODE {\em(\ref{ode})} is trivialisable by contact 
transformations if and only
if the invariants {\em(\ref{doubrov1})} 
vanish, and the following conditions hold
\begin{itemize}
\item[$n=4$.] $F_{333}=6F_{233}+F_{33}^2=0$,
\item[$n=5$.] $F_{44}=6F_{234}-4F_{333}-3F_{34}^2=0$,
\item[$n=6$.] $F_{55}=F_{45}=0$,
\item[$n\geq 7$.] $F_{(n-1)(n-1)}=F_{(n-1)(n-2)}=F_{(n-2)(n-2)}=0$.
\end{itemize}
\end{theo}
\section{Twistor theory}
\setcounter{equation}{0}
\label{sec_twistor}
We shall now show how, in the real analytic case, 
the paraconformal structure (\ref{paracon})
on $M$ can be encoded in a holomorphic geometry of rational curves 
embedded in a complex surface ${\cal Z}$. The $n$th order ODE on 
${\cal Z}$ will then implicitly be given by the embedding 
$L\subset {\cal Z}$, provided that $L$ has self-intersection number $n$.

In this section we regard (\ref{relation}) as a holomorphic 
relation between complex coordinates $(x, y, {t})$, 
and examine the geometry of the complexified  correspondence space 
(\ref{correspondence}) and the associated double fibration picture.
The surface ${\cal Z}$ will eventually become a twistor space of
(\ref{ode}). We shall  however set up a more general correspondence, and
consider Legendrian curves in three-dimensional
twistor spaces.

Let $Y$ be a complex three-fold with an embedded rational curve $L$
with a normal bundle $N=\O(n-2)\oplus\O(n-2)$. We have
$H^1(\CP^1,  \O(n-2)\oplus\O(n-2))=0$, and so
the moduli space of such curves in $Y$ is a manifold ${\cal M}$ of 
dimension equal to 
\[
\dim H^0(\CP^1,  \O(n-2)\oplus\O(n-2)))=
2n-2.
\] 

Now we restrict our attention to a moduli space $M$ of contact 
(Legendrian) curves with normal bundle $N$. The canonical line bundle
of holomorphic three-forms on $Y$ restricted to a curve $L$ is 
\[
\kappa(Y)=T^*L\otimes \Lambda^2(N^*)=\O(2-2n),
\]
since $T^*\CP^1=\O(-2)$.
From the general theory of contact structures it follows that
the contact line bundle is given by ${L_c}^2=\kappa(Y)$.
Now pick a section  of ${L_c}^*$ (a contact one-form), and contract it
with a tangent vector  to a rational curve to get a section of
\[
(\O(1-n)\otimes\O(2))^*.
\]
The vanishing of this section (the Legendrian condition) gives
$\dim H^0(\CP^1,  \O(n-3))=n-2$ conditions on ${\cal M}$. Therefore
the dimension of the moduli space $M$ of Legendrian curves is
\[
\dim M= (2n-2)-(n-2)=n.
\]
This can be summarised by the double fibration picture
\be
\label{doublefib2}
{M}\stackrel{\hat{p}}\longleftarrow 
\hat{\cal{F}}\stackrel{\hat{q}}\longrightarrow {Y}.
\ee
The curves $\hat{q}(\hat{p}^{-1}({\bf t}))\cong \CP^1$ are Legendrian with
respect to the contact form on $Y$. 
Bryant's generalisation \cite{B91} of the Kodaira theorems
guarantees that the moduli space $M$ of Legendrian rational curves is stable
under small deformations of $Y$.

 Consider the special case $Y=\PP(T\z)$. A rational curve $L$ with normal bundle 
$\O(n-1)$ in $\z$ has a natural lift $\hat{L}$ to $Y$, given by
$z\in L\rightarrow (z, \dot{z}\in T_zL)$. The lifted curves are
Legendrian with respect to the canonical contact structure on the
projectivised tangent bundle.
The double fibration (\ref{doublefib2}) reduces to (\ref{doublefib}).

The existence of the complexified paraconformal structure (\ref{paracon})
follows 
from the structure of the normal bundle. From Kodaira \cite{Ko63} theory,
since the appropriate obstruction groups vanish, we have
\begin{equation}
\label{solder}
T_tM=\Gamma(L_t, N_t)=S^{n-1}(\spp_t), \qquad \spp_t=\C^2,
\end{equation}
where $N_t$ is the normal bundle to the rational curve $L_t=\CP^1$ in
${\cal Z}$ corresponding to the point $t\in M$.  
The nontrivial examples of ODEs satisfying all the constraints
(\ref{invariantsC}) can therefore be constructed by
applying algebraic operations on a rational curve embedded in
a total space of ${\cal O}(N)$ for $N$ sufficiently large \cite{DGM05}.

The correspondence space 
${\cal F}=M\times\CP^1$ is equipped with a canonical $(n-1)$
dimensional distribution ${\cal D}$, such that ${\cal Z}={\cal F}/{\cal D}$.
The normal bundle to a rational curve $L_t:=q(p^{-1}({\bf t}))$ 
consists of vectors tangent to $M$ at ${t}$ (horizontally lifted to
$T_{ {t}, \lambda} {\cal F}$) modulo ${\cal D}$. 
Therefore we have a sequence of sheaves over $\CP^1$
\[
0\longrightarrow {\cal D} \longrightarrow \C^n \longrightarrow
{\cal O}(n-1)\longrightarrow 0.
\]
The map $\C^n \longrightarrow \O(n-1)$ is given by
$V^{A_1A_2...A_{n-1}}\longrightarrow 
V^{A_1A_2...A_{n-1}}z_{A_1}z_{A_2}...z_{A_{n-1}}$.  
Its kernel consists of
vectors of the form $z^{(A_1}\lambda^{A_2...A_{n-1})}$ with 
$\lambda^{A_2...A_{n-1}}\in \C^{n-1}$  varying. The
twistor distribution is therefore ${\cal D}=\O(-1)\otimes
\mbox{S}^{(n-2)}(\C^2)$. This distribution is the geodesic spray
(\ref{gspray}) if $n=2$, 
or the Einstein--Weyl Lax pair (\ref{EWlax}) if $n=3$.
\section{Exotic ${\cal G}_3$ holonomy and fourth order ODEs}
\setcounter{equation}{0}
\label{sec_holo}
In this Section we shall make contact with Bryant's work \cite{B91}
and show that if $n=4$ the W\"unschmann conditions (\ref{invariantsC}) 
are equivalent to
the existence of certain exotic holonomy on $M$.

Recall the notation introduced at the beginning of Section (\ref{sec_wun}) 
and define ${\cal G}_k\subset GL(V_k)$ by
\[
{\cal G}_{k}=\{g\in GL(V_k)|\; g({\bf t})\in C_k \;\;\mbox{if}\;\;{\bf
t}\in C_{k}\}.
\]
Can ${\cal G}_{n-1}$ appear as a holonomy group of a torsion-free
connection of an $n$--dimensional manifold? 
Bryant \cite{B91} 
has examined  Berger's criteria, and established that 
the answer is `no' if $n>5$ (the case $n=5$ is special, as
the five-dimensional representation of $SL(2, \C)$ is the holonomy of a
symmetric space $M=SL(3, \C)/SL(2, \C)$). 

Let  $<, >_2$  be given by (\ref{inner_product}), and let
\[
g(X, X, X, X):=<<X,X>_2, <X, X>_2>_2
\]
be a ${\cal G}_3$ invariant quartic form on $TM$.

\begin{defi}
The vector $X\in T_tM$ is null iff
\[
Q(X):=g(X, X, X, X)=0.
\]
The $\a$--plane is a two--dimensional plane in $T_tM$ spanned by
vectors $X, Y$ such that
\[
Q(X+\l Y)=0 
\]
for each value of a parameter $\l$.
\end{defi}
\label{a_planes}
The null vectors in this sense 
correspond to polynomials of the form ${\bf p}^2 {\bf
 r}$. Vectors in an $\a$--plane are then obtained by varying ${\bf r}$
 and keeping ${\bf p}$ fixed.
\begin{theo}[Bryant \cite{B91}]
\label{theo_br}
A four--dimensional manifold $M$ admits a torsion free connection with
holonomy ${\cal G}_3$ iff for every $\a$--plane there exists a
two--dimensional
surface $\Sm\subset M$ (called $\a$--surface) tangent to this
$\a$--plane.

The space of torsion--free ${\cal G}_3$ structures modulo
diffeomorphisms depends on four arbitrary functions of three variables.
\end{theo}
Bryant has also shown \cite{B91} that if $M$ admits a ${\cal G}_3$ structure, 
then there exists a three--parameter family of $\a$--surfaces.
In the complexified category this family is parametrised by points of
a complex three-fold ${Y}$, and the Legendrian $\O(2)\oplus\O(2)$ 
curves in ${Y}$ correspond to points in $M$ (compare this
with the twistorial treatment in Section \ref{sec_twistor}).

{\bf Proof of Theorem \ref{theo_G3}.}
The conditions (\ref{W4_1}, \ref{W4_2}) are the invariants
(\ref{invariantsC}) with $n=4$ which establishes the equivalence of
$(1)$ and $(3)$.

To show $(1)\rightarrow (2)$ observe that the existence of the
paraconformal structure (\ref{paracon}) implies the existence of a
symplectic structure $\varepsilon$ up to scale on each spin space.
This gives us the symmetric quartic form 
\[
g(X, X, X, X)=
\varepsilon^{A_2B_2}\varepsilon^{A_3B_3}
\varepsilon^{C_2D_2}\varepsilon^{C_3D_3}
\varepsilon^{A_1C_1}\varepsilon^{B_1D_1}
X_{A_1A_2A_3}X_{B_1B_2B_3}X_{C_1C_2C_3}X_{D_1D_2D_3}.
\]
The quartic $Q(X)=0$ selects null vectors, and $\a$--planes.
Theorem \ref{theo_br} asserts that these planes are 
integrable iff they come from a ${\cal G}_3$ structure. 
But they will always be integrable in the paraconformal case because 
of the following interpretation:
Fixing a point in $Z\in{\cal Z}$ gives a three--dimensional surface 
in $N\subset M$, such that its normal
$\nabla Z$ is a perfect cube. Fixing a point and a direction in ${\cal Z}$
gives an $\a$--surface $\Sm\subset N$ (think of $Z, Z', Z'', Z'''$ 
as coordinates in $M$). This corresponds to fixing a point in 
${Y}=\PP(T^*{\cal Z})$, and to Bryant's $\a$--plane with normals
$Z, Z'$ which are gradients (so that it really is a surface).

It remains to demonstrate $(2)\rightarrow(1)$.
Suppose we look directly for the quartic as 
\be
\label{quartic}
g=g_0+g_1, 
\ee
where $g_0$ is the form when the ODE is trivial $(F=0)$ and $g_1$ 
is a combination of all possible
lower order terms (lower order in the sense of fewer derivatives of $y$):
\begin{eqnarray*}
g_0&=&18\d y\;\d p\;\d q\;\d r-9 (\d y)^2\;(\d r)^2+3(\d p)^2\;(\d q)^2
-8(\d p)^3\;\d r- 6\d y\;(\d q)^3,\\
g_1&=&\a (\d y)^4+\beta (\d y)^3\;\d p+\gamma (\d y)^2\;(\d
p)^2+\delta (\d y)^3\;(\d q)+\epsilon \d y\; (\d p)^3\\
&&
+\xi (\d y)^2\; \d p\;\d q +\eta (\d y)^3 \;\d r+\kappa (\d
p)^4
+\gamma (\d y) \;(\d p)^2\; \d q+\mu (\d y)^2 (\d q)^2+\nu
(\d y)^2 \;\d p\;\d r \\
&&+\zeta (\d p)^3 \;\d q+\pi \d y \;\d p\; (\d
q)^2 
+\theta \d y\;(\d p)^2\;\d r+\phi (\d y)^2\;\d q\;\d r,
\end{eqnarray*}
and $p=y', q=y'', r=y''', s= y''''=F(x, y, p, q, r)$.

We need to fix 15
coefficients $(\a, ..., \phi)$. We impose 
\[
g'=\Lambda g
\] 
for some $\Lambda$ and work systematically through
the coefficients, fixing them in order (see Appendix for the details
of this calculation).  We solve equations (\ref{eeq5}--\ref{eeq14})
and (\ref{eeq15}, \ref{eeq17}). Now equation (\ref{eeq16}) 
becomes (\ref{W4_2}).
Then we solve (\ref{eeq18}, \ref{eeq20}). Now (\ref{eeq19}) and
(\ref{W4_2}) give (\ref{W4_1}).

The remaining conditions (\ref{eeq21}--\ref{eeq24}) give two more conditions
on $F$, but we know these will be satisfied because of the
paraconformal argument: once we impose (\ref{W4_1},\ref{W4_2})
the quartic will exist
and we can check that the coefficients we have found by a direct approach agree
with what the paraconformal method has told us. Thus (\ref{W4_1},\ref{W4_2}) 
are necessary and sufficient for integrability. 
\koniec
Readers who compare our treatment with that of Bryant's  will recall 
(Theorem 4.5 in \cite{B91})) that the existence
of a torsion-free ${\cal G}_3$-structure on the moduli space was equivalent
to the vanishing of two primary invariants $(A, C)$, and two secondary
invariants  $(B, D)$ which are well defined only if $A=C=0$.
All these invariants  are polynomials in the third order derivatives 
of the function $F$. Bryant has pointed out 
\cite{Bprivate} that  two conditions 
(\ref{W4_1}, \ref{W4_2}) are equivalent to the vanishing of A and B. 
If $A=B=0$ then  C and D vanish identically.

\section{Torsion--free paraconformal connections}
\setcounter{equation}{0}
\label{sec_conn}
Let $M$ admit a paraconformal structure, and let
\[
D:\Gamma(\spp^k)\longrightarrow \Gamma(\spp^k\otimes T^*M)=
\Gamma(\spp^{k+n-1}),\qquad \mbox{where}\;\;\spp^k=\spp^{\otimes k},\;\;k=0, 1, ..., n-1
\]
be a connection. In this section we shall show that if $D$ preserves
the paraconformal structure (\ref{paracon}) then it necessarily has torsion
if $n=4$ or $n\geq 6$.
\begin{defi} 
The connection $D$ is called paraconformal if
its action  on elements of $V_1$ is given by
\[
D{\bf p}=U\otimes{\bf p}+V\otimes{\bf q}
\]
for some $U, V\in\spp^{n-1}$. 
\end{defi}
This definition implies that
$
D{\bf E}=(n-1)U\otimes{\bf E}+V\otimes{\bf E}'.
$
The torsion free condition becomes
\[
U=A{\bf E}+B{\bf E}',\;V=(n-1)B{\bf E}+C{\bf E}'
\]
for some $A, B, C\in V_0$. 
Demanding that the connection does not depend
on $(x, y)$ gives  the consistency condition
\be 
\label{diff_constr}
(D{\bf p})''=D({\bf p}''),
\ee 
which yields
\begin{eqnarray*}
0&=&(D{\bf p})''-D(P{\bf p}+Q{\bf q})=(\frac{\d^2}{\d
  x^2}-Q\frac{\d}{\d x}-P)D{\bf p}-(DP)\otimes{\bf p}-(DQ)\otimes{\bf
  q} \\
&=&(\a_1{\bf E}+ \a_2{\bf E}'+\a_3{\bf E}''+\a_4{\bf E}'''
)\otimes {\bf p}+
(\beta_1{\bf E}+ \beta_2{\bf E}'+\beta_3{\bf E}''+\beta_4{\bf E}''')
\otimes {\bf q}\\
&&-
(\sum_{i=0}^{n-1}{\bf E}^{(i)}\frac{\p P}{\p y^{(i)}}\otimes{\bf p}+
\sum_{i=0}^{n-1}{\bf E}^{(i)}\frac{\p Q}{\p y^{(i)}}\otimes{\bf q}).
\end{eqnarray*}
The coefficients of ${\bf E}^{(i)}\otimes {\bf p}$, and
${\bf E}^{(i)}\otimes {\bf q}$ have to vanish, which gives 
$2n$ conditions on $F$.
Here  $(P, Q)$ are given by (\ref{pandq}) and
$\a_1, ..., \a_4, \beta_1, ..., \beta_4$ can be determined in terms of
$A, B, C, P, Q$, and their derivatives. \vskip 10pt\noindent
{$\bf n=2.$}  In this case ${\bf E''}$ and
${\bf E'''}$ are determined in terms of ${\bf E}$ and ${\bf E'}$ according 
to (\ref{expE}), and therefore $\a_3, \a_4, \beta_3, \beta_4$ all vanish.
There are no conditions on  $F$  arising from (\ref{diff_constr}).
However imposing the geodesic conditions on curves $Z=const$ yields
(\ref{geodesiccond}). Let $t^A$ and $z^A=\dot{t}^A$ be 
local coordinates on $M$ and $T_tM$. The Christoffel
symbols $\Gamma_{AB}^C=\Gamma_{AB}^C({t})$ of $D$ are defined up
to  projective equivalence 
\[
\Gamma_{AB}^C\sim \Gamma_{AB}^C+\delta_{(A}^C\omega_{B)}
\]
for some $\omega_B=\omega_B({t})$. Let $t^C=t^C(\tau)$ be solutions
to 
\be
\label{geodeq}
\ddot{t}^C+\Gamma_{AB}^C\dot{t}^A\dot{t}^B=v\dot{t}^C, 
\qquad {\dot{}}=\frac{\d}{\d\tau}.
\ee
where $v$ is some function.
These geodesic curves lift to the integral curves of the geodesic
spray which is a projection of
\be
\label{gspray}
L=z^A\frac{\p}{\p t^A}-\Gamma_{AB}^Cz^Az^B\frac{\p}{\p z^C}
\ee
from $TM$ to $\PP(TM)$.
Eliminating $\tau$ from (\ref{geodeq}) leads to a second order ODE for
$t^1=t^1(t^0)$ which is at most cubic in the first derivatives (the
cubic term is given by $\varepsilon_{CD}\Gamma_{AB}^Cz^Az^Bz^D$, where 
$z^A$ are homogeneous coordinates on $\PP(TM)$, and $z^1/z^0=\d t^1/
\d t^0$). This ODE is dual to (\ref{ode}) in the sense of Cartan \cite{C22}. 
It could also be read off from the relation (\ref{relation})
by rewriting it as $t^1=K(t^0, x, y)$, and eliminating $(x,y)$
between $K$ and its first two derivatives w.r.t $t^0.$
A second order ODE (\ref{ode}) is trivialisable by point transformations iff
the curvature of the projective connection vanishes. This curvature vanishes
iff (\ref{geodesiccond}) holds, and $F_{1111}=0$.
\vskip 10pt\noindent
{$\bf n=3$}. Now $\a_4$ and $\beta_4$ vanish.
The compatibility conditions in (\ref{diff_constr})
fix $A, B, C$. Imposing the totally geodesic condition $C=0$
on the null surfaces $Z=const$  gives 
the constraint (\ref{EWcond}). We shall come back to this case in the next
Section.\vskip 10pt\noindent
${\bf n\geq 4.}$
The coefficients of 
${\bf E}'''\otimes{\bf p}, {\bf E}'''\otimes{\bf q}$ 
and ${\bf E}''\otimes{\bf p}$
fix $(A, B, C)$ in terms of $F$ and its derivatives and 
the coefficients of ${\bf E}''\otimes{\bf q}, {\bf E}'$ and ${\bf E}$ 
give five equations (\ref{A1}--\ref{A5})  for $F$ (the details are
in the Appendix). 
In particular (\ref{A3}) yields
\be
\label{3rdcond}
\Big(\frac{6n-8}{n}\Big)F_{n-1}F_{(n-1)3}+(8-2n)\frac{\d}{\d x}F_{(n-1)3}
-(2n-2)F_{(n-1)2}+6F_{(n-2)3}=0.
\ee
{\bf Proof of Theorem \ref{connection}}.
If $n=4$ the condition (\ref{3rdcond}) reduces to
\[
F_{33}=0,
\]
and the W{\"u}nschmann conditions (\ref{W4_1}, \ref{W4_2}) now imply
$A=B=C=0$. Some computer algebra reduces (\ref{W4_1}, \ref{W4_2}) to
\[
F=\a(x)+\beta(x)y+\gamma(x)y'+\delta(x)y''+\epsilon(x)y''',
\]
where $\a(x), \delta(x), \epsilon(x)$ are arbitrary, and 
\[
\beta=\frac{11}{1600}\epsilon^4-\frac{9}{50}\epsilon^2\epsilon'-\frac{1}{200}\epsilon^2\delta+\frac{21}{100}(\epsilon')^2+\frac{1}{50}\epsilon'\delta
-\frac{9}{100}\delta^2+\frac{7}{20}\epsilon\delta''-\frac{1}{5}
\epsilon'''+\frac{3}{10}\delta''-\frac{1}{4}\epsilon\delta'
\]
\[
\gamma=\frac{3}{4}\epsilon\epsilon'-\frac{1}{2}\epsilon''+\delta'
-\frac{1}{2}\delta\epsilon-\frac{1}{8}\epsilon^3.
\]
We can however preform a point transformation (which is in fact
fibre preserving)
\[
y=a(x)\hat{y}(x)+b(x),\qquad x=c(\hat{x}),
\]
and choose the functions $(a, b, c)$ to set $\a=\delta=\epsilon=0$.
The resulting 4th order ODE (\ref{ode}) is therefore trivial up to point 
transformations.

If $n>4$, then  
the coefficients
of ${\bf E}^{(k)}$ with $k>3$ give
\be
\label{other_cond}
\frac{\p P}{\p y^{(k)}}=\frac{\p Q}{\p y^{(k)}}=0,\qquad k=4, 5, ..., n-1.
\ee
The conditions  (\ref{other_cond}) give the following
\begin{enumerate}
\item $n=6$
\be
\label{eq45}
F_{55}=F_{45}=3F_{44}-8F_{53}=0,
\ee
\item $n\geq 7$
\be
\label{eq46}
F_{(n-1)k}=F_{(n-2)(k+1)}=3F_{(n-2)4}-(n-2)F_{(n-1)3}=0,\qquad
k=4, ..., n-1.
\ee
\end{enumerate}
Equations (\ref{invariantsC}) and (\ref{A1}--\ref{A5}) 
also have to be satisfied.  The equations 
(\ref{eq45}, \ref{eq46}) and Theorem \ref{doubrov2}  imply that
the ODE (\ref{ode}) is trivialisable by contact transformation if $n\geq 6$. 
\koniec
The case $n=5$ is exceptional. The conditions (\ref{other_cond})
give
\be
\label{eq44}
F_{44}=0,
\ee
but this is not sufficient to guarantee the trivialisability.
In fact the five parameter family of conics in
the complex projective plane gives a counterexample.
In this case $M$  (a real form of $PSL(3,\C)/SO(3,\C)$) 
is a 5-dimensional space with 
the paraconformal structure
which admits a torsion free paraconformal connection.
 The 5-parameter family of conics in $\CP^2$ is given by (\ref{relation})
with
\[
\Psi=y^2+t^1x^2+2t^2xy+t^3x+t^4y-t^5,
\]
where $(y, x)$ are inhomogeneous coordinates on $\CP^2$.
We regard $y$ as a function of $x$, and implicitly 
differentiate the relation (\ref{relation}) five times w.r.t $x$.
Solving for $y^{(5)}$ yields the desired $5$th order ODE
\be
\label{order_5}
y^{(5)}=-\frac{40 r^3}{9 q^2}+5\frac{rs}{q} ,
\ee
where  $y'=p, y''=q, y'''=r, y''''=s$. 

The W\"unschmann conditions
(\ref{W4_1}, \ref{W4_2}), as well as the 
six conditions (\ref{other_cond}, \ref{A1}, \ref{A2}, \ref{A3}, \ref{A4}, 
\ref{A5} ) hold, 
so the paraconformal torsion
free connection  exists in this case. 
The equation (\ref{order_5}) 
is nevertheless not contact equivalent to $y^{(5)}=0$, as
the  invariant $6F_{234}-4F_{333}-3F_{34}^2$ 
from Theorem \ref{doubrov2} doesn't vanish, and is equal to 
$(5/3)q^{-2}$.
\section{From Einstein--Weyl structures to third order ODEs}
\label{sec_EW}
\setcounter{equation}{0}
In three dimensions the existence a 
paraconformal structure (\ref{paracon}) 
is equivalent to the existence of a conformal
structure $[h]$ of signature $(++-)$. This is a well known fact based
on representing vectors as symmetric matrices
\[
X^a=(X^1, X^2, X^3)\longrightarrow X^{AB}=\left (
\begin{array}{cc}
X^1+X^2&X^3\\
X^3&X^1-X^2
\end{array}
\right )\in \Gamma (\spp\otimes\spp),
\]
where $X^a$ are components of $X$ w.r.t. some basis.

The matrices corresponding to null vectors (i.e $h(X, X)=0, h\in [h]$) 
have vanishing determinant,
and must have rank one. Therefore $X^{AB}=p^Ap^B$ for such vectors.

Set $n=3$, and assume that the 3rd order ODE (\ref{ode}) satisfies the 
W\"unschmann condition (\ref{wun3}).
For each choice of $(x, y)$ (\ref{relation})
defines a surface in $M$ which is null w.r.t $[h]$.
In the last section we have shown that if $n=3$
the null surfaces $y=Z(x, t^a)$ 
are totally geodesic w.r.t some torsion--free
connection $D$ if $F(x, y, y', y'')$ satisfies the constraint
(\ref{EWcond}), 
and it is well known \cite{C43, H82, PT93} 
that the Einstein--Weyl (EW) equations are equivalent to 
the existence of a two  dimensional family of
surfaces $\Sm\subset M$ which are 
null with respect to $[h]$, and  totally geodesic with respect to $D$.

Let $M$ be a  $3$-dimensional manifold with a torsion-free connection $D$,
and a conformal structure $[h]$ of signature $(++-)$ 
which is compatible with $D$ in a sense
$Dh=\om\otimes h $
for some one-form $\om$.
Here $h\in[h]$ is a representative metric in a conformal class. If 
we  change this representative by $h\rightarrow \psi^2 h$, 
then $\om\rightarrow \om +2\d \ln{\psi}$, where $\psi$ is a 
non-vanishing function on $W$. A triple $(M, [h], D)$ is called a Weyl 
structure.
The conformally invariant Einstein--Weyl (EW) equations state
that the symmetrised part of the Ricci tensor of $D$ is proportional
to the representative of $[h]$.

Given a third order ODE which satisfies (\ref{wun3}, \ref{EWcond}),
 the EW structure can be reconstructed following
the steps described in \cite{C43, T00}. The problem of reconstructing the ODE
starting from a given EW structure was left open in these
references. We shall present a method which reduces the problem of finding
the allowed ODE to a system of linear PDEs.
First recall the Lax representation
for the EW equations \cite{DMT00}. 
Let $X_1, X_2, X_3$ be three independent 
vector fields on $M$, and let $e_1, e_2, e_3$ be the dual one-forms.
Assume that 
\[
h=e_2\otimes e_2-2(e_1\otimes e_3+e_3\otimes e_1)
\] 
and some one--form $\om$ together give an EW structure.
Let $X(\ll)=X_1-2\ll X_2+\ll^2X_3$ where $\ll\in\CP^1$ is a 
projective coordinate on the fibres of $\spp\rightarrow M$. 
Then $h(X(\ll), X(\ll))=0$ for 
all $\ll\in\CP^1$ so $X(\ll)$ determines a sphere of null vectors. 
The vectors $X_1-\ll X_2$ and
$X_2-\ll X_3$ form  a basis of the orthogonal complement
of $X(\ll)$. For each $\ll\in\CP^1$ they span 
a null two-surface. 
Therefore  the Frobenius theorem implies that the 
horizontal lifts to $\spp$
\be
\label{EWlax}
L_0=X_1-\ll X_2+l_0\p_{\ll},\qquad
L_1=X_2-\ll X_3+l_1\p_{\ll}
\ee
satisfy $[L_0, L_1]=\a L_0+\beta L_1$
for some $\a, \beta$ which are linear in $\ll$.
The functions $l_0$ and  $l_1$ 
are third order in $\ll$, because the M{\"o}bius transformations of
$\CP^1$ are generated by  vector fields quadratic in $\ll$.

To find the third order ODE corresponding to $([h], D)$ we construct
two independent solutions $x(t_a, \l), y(t_a, \l)$ to the pair of linear PDEs
\[
L_0f=L_1f=0,
\]
and eliminate $\ll$ between $x$ and $y$. This gives $y=Z(x, t^a)$. Now
we follow the prescription given in the introduction to produce 
the third order ODE. Both invariants (\ref{wun3}) and
(\ref{EWcond}) will be satisfied as a consequence of the EW condition.

As an example, 
consider the Einstein--Weyl $(++-)$ structure  on Thurston's Nil manifold 
$S^1\times \R^2$ \cite{PT93, T00} given by 
\[
h=\a^2(\d t^2+t^1\d t^3)^2-4\d t^1\d t^3,\qquad \om=\a^2(\d t^2+ t^1\d t^3).
\]
Choose the Lax pair 
\[
L_0=\p_1+\a^{-1}\ll\p_2, \qquad L_1=-\a^{-1}\p_2-\ll(\p_3-t^1\p_2)+\a\l\p_\l
\]
so that $[L_0, L_1]=0$.
We find a kernel of $(L_0, L_1)$  to be 
\[
x=\ll+\a t^3,\qquad  y=\ll t^1-\a t^2 -\a^{-1}\ln{\ll}
\]
so that the totally geodesic surfaces are given by $y=Z(x, t^a)$ with 
\[
Z(x, t_a)=(x-\a t^3)t^1-\a t^2-\a^{-1}\ln{(x-\a t^3)}.
\]
The resulting third order ODE is
\[
y'''=2\sqrt{\a}(y'')^{3/2}
\]
which is a special case of our general example (\ref{example1}) with $n=3$,
and $(t^1, t^2, t^3)$ redefined.

\section*{Acknowledgements}
We thank Robert Bryant and Boris Doubrov for useful discussions
and correspondence which resulted in many improvements. We also thank
Michael Eastwood for bringing reference \cite{Dou} to our attention, and 
the anonymous referee for valuable comments.
\section*{Appendix}
\setcounter{equation}{0}
\appendix
\def\theequation{\thesection{A}\arabic{equation}}
{\bf Determining P and  Q.} 
The first step is to calculate recursive formulae for $a_{ij}$. This 
yields $a_{00}=1$, and
\begin{eqnarray*}
a_{ik}&=&0,\qquad \mbox{for}\qquad k>i,\\
a_{(i+1) k}&=&(a_{ik})'+kQa_{ik}+(n-k)a_{i(k-1)}+(k+1)a_{i(k+1)}P,\qquad
0\leq k\leq i-1,\\
a_{(i+1)i}&=&(a_{ii})'+iQa_{ii}+(n-i)a_{i(i-1)},\\
a_{(i+1)(i+1)}&=&(n-i-1)a_{ii}.
\end{eqnarray*}
These relations give
\[
a_{ii}=\frac{(n-1)!}{(n-i-1)!}, \qquad a_{(i+1)i}=\frac{(i+1)i}{2}Qa_{ii},
\qquad a_{n(n-2)}=\alpha P+\beta Q^2+\gamma Q',
\]
where
\[
\alpha=\frac{n(n+1)!}{6},\qquad \beta=\frac{n!}{24}(3n-5)(n-1)(n-2),\qquad
\gamma=\frac{n!}{6}(n-1)(n-2).
\]
Solving the last two equations in (\ref{matrixODE}) corresponding to 
$j=n-1$, and $j=n$ for $P$ and $Q$ yields (\ref{pandq}). \vskip 10pt\noindent
{\bf Conditions for F.} We shall give explicit 
forms of  (\ref{invariantsC}) for $n=5$ 
\begin{eqnarray*}
0&=&-\frac{1}{5}\,{\frac {\d^{4}}{\d{x}^{4}}}F_{{4}}-{\frac {2}{25}}\,F_{{3}}{
\frac {\d^{2}}{\d{x}^{2}}}F_{{4}}-{\frac {7}{25}}\,\left ({\frac {\d}{
\d x}}F_{{4}}\right ){\frac {\d}{\d x}}F_{{3}}-{\frac {28}{3125}}\,
\left (F_{{4}}\right )^{5}+{\frac {16}{25}}\,\left ({\frac {\d}{\d x}}
F_{{4}}\right ){\frac {\d^{2}}{\d{x}^{2}}}F_{{4}}\\
&& -F_{{0}}+
{\frac {8}{25}}\,F_{{4}}{\frac {\d^{3}}{\d{x}^{3}}}F_{{4}}-\frac{1}{5}\,F_{
{2}}F_{{3}}-{\frac {7}{100}}\,F_{{2}}\left (F_{{4}}\right 
)^{2}+\frac{1}{5}\,F_{{2}}{\frac {\d}{\d x}}F_{{4}}-{\frac {11}{125}}\,F_{{
4}}\left (F_{{3}}\right )^{2}\\
&&-{\frac {141}{2500}}\,F_{{3}}
\left (F_{{4}}\right )^{3}+
{\frac {137}{1250}}\,\left (F_{{4}}
\right )^{3}{\frac {\d}{\d x}}F_{{4}}-{\frac {9}{25}}\,F_{{4}}
\left ({\frac {\d}{\d x}}F_{{4}}\right )^{2}-{\frac {9}{50}}\,F_{{4}}
{\frac {\d^{2}}{\d{x}^{2}}}F_{{3}}\\
&&
-{\frac {103}{500}}\,\left (F_{{4
}}\right )^{2}{\frac {\d^{2}}{\d{x}^{2}}}F_{{4}}+{\frac {101}{1000
}}\,\left (F_{{4}}\right )^{2}{\frac {\d}{ \d x}}F_{{3}}+{\frac {7}{
50}}\,F_{{3}}{\frac {\d}{\d x}}F_{{3}}+\frac{1}{5}\,{\frac {\d^{3}}{\d{x}^{3}
}}F_{{3}}+{\frac {28}{125}}\,F_{{4}}F_{{3}}{\frac {\d}{\d x}}F_{
{4}}\\
0&=&
-12\,F_{{2}}-{\frac {36}{5}}\,F_{{3}}F_{{4}}-{\frac {48}{25}}
\,\left (F_{{4}}\right )^{3}+{\frac {72}{5}}\,F_{{4}}
{\frac {\d}{
\d x}}F_{{4}}-12\,{\frac {\d^{2}}{\d{x}^{2}}}F_{{4}}+18\,{\frac {\d}{
\d x}}F_{{3}}\\
0&=&
{\frac {102}{25}}\,F_{{4}}{\frac {\d^{2}}{\d{x}^{2}}}F_{{4}}+{
\frac {18}{5}}\,{\frac {\d^{2}}{\d{x}^{2}}}F_{{3}}-{\frac {16}{5}}\,{
\frac {\d^{3}}{\d{x}^{3}}}F_{{4}}+{\frac {68}{25}}\,\left ({\frac {\d}
{\d x}}F_{{4}}\right )^{2}\\&& -{\frac {4}{625}}\left (F_{{4}}\right 
)^{4}-{\frac {16}{25}}\,\left (F_{{3}}\right )^{2}-{\frac {34}{125}
}\,F_{{3}}\left (F_{{4}}\right )^{2}+{\frac {4}{25}}\,F_{{3}}
{\frac {\d}{\d x}}F_{{4}}-{\frac {172}{125}}\,\left (F_{{4}}\right 
)^{2}{\frac {\d}{\d x}}F_{{4}}\\
&&-4\,F_{{1}}-\frac{2}{5}\,F_{{2}}F_{{4}}
-\frac{9}{5}\,F_{{4}}{\frac {\d}{\d x}}F_{{3}}=0.
\end{eqnarray*}
The conditions for $n>5$ can be written down using
recursive relations and MAPLE, but the resulting formulae are very
long, and inconclusive. \vskip 10pt\noindent
{\bf Calculations leading to a proof of Theorem \ref{theo_G3}}.
We collect the terms in $g'-\Lambda g=0$ by order of derivatives
(eg $\d y\;\d p\;\d r$ has order $4=0+1+3$). \\
6th order
\be\label{eeq5}
18 F_3+2\pi+2\theta+2\phi-18\Lambda=0
\ee
\be\label{eeq6}
-18F_3+\phi+9\Lambda=0
\ee
\be \label{eeq7}
\pi+6\Lambda=0
\ee
\be\label{eeq8}
-8F_3+\zeta+\theta+8\Lambda=0
\ee
\be\label{eeq9}
3\zeta+\pi-3\Lambda=0.
\ee
5th order
\be\label{eeq10}
-18F_2+2\mu+\nu+\phi'+\phi F_3-\Lambda\phi=0
\ee
\be \label{eeq11}
18F_2+2\l+2\mu+\pi'-\Lambda\pi=0
\ee
\be\label{eeq12}
\l+2\nu+\theta F_3+\theta'-\Lambda\theta=0
\ee
\be\label{eeq13}
-8F_2+4\kappa+\l+\zeta' -\Lambda\zeta=0.
\ee
4th order
\be\label{eeq14}
-18 F_1+\xi+3\eta+\nu F_3+\nu'-\Lambda\nu=0
\ee
\be \label{eeq15}
\xi+\phi F_2+\mu'-\Lambda\mu=0
\ee
\be\label{eeq16}
18 F_1+3\epsilon+2\xi+\theta F_2+\l'-\Lambda\l=0
\ee
\be\label{eeq17}
-8F_1+\epsilon+\kappa'-\Lambda\kappa=0.
\ee
3rd order
\be\label{eeq18}
-18F_0+\delta+\eta F_3+\eta'-\Lambda\eta=0
\ee
\be\label{eeq19}
18 F_0+2\gamma+3\delta+\nu F_2+\phi F_1+\xi'-\Lambda\xi=0
\ee
\be\label{eeq20}
-8F_0+2\gamma+\theta F_1+\epsilon'-\Lambda\epsilon=0.
\ee
2nd order
\be\label{eeq21}
3\beta+\nu F_1+\theta F_0+\gamma'-\Lambda\gamma=0
\ee
\be\label{eeq22}
\beta+\eta F_2+\phi F_0+\delta'-\Lambda\delta=0.
\ee
1st order
\be\label{eeq23}
4\a +\eta F_1+\nu F_0+\beta' -\Lambda\beta=0.
\ee
0th order
\be
\label{eeq24}
\eta F_0+\a'-\Lambda\a=0.
\ee
\vskip 10pt\noindent
{\bf Conditions for F leading to a proof of Theorem \ref{connection} .}
\[
B=P_3,\qquad C=Q_3,\qquad 
A=P_2+QP_3-2Q_3P-2\frac{\d}{\d x}P_3,
\]
where $P, Q$ are given by (\ref{pandq}).
\be
\label{A1}
2\frac{\d A}{\d x}-Q_0+(n-1)\frac{\d}{\d x}(BQ+\frac{\d}{\d x}B)=0
\ee
\be\label{A2}
\frac{\d^2}{\d x^2}C+2A-Q_1+\frac{\d}{\d x}(CQ)+(n-1)QB
+2n\frac{\d B}{\d x}=0
\ee
\be
\label{A3}
CQ+2\frac{\d C}{\d x}-Q_2+(n+1)B=0
\ee
\be
\label{A4}
\frac{\d^2 A}{\d x^2}-Q\frac{\d A}{\d x}-P_0+(n-1)\frac{\d}{\d x}(BP)
+(n-1)P\frac{\d B}{\d x}=0
\ee
\be
\label{A5}
2\frac{\d A}{\d x}-P_1-QA+\frac{\d^2 B}{\d x^2}-Q\frac{\d B}{\d x}+
\frac{\d}{\d x}(CP)+P\frac{\d C}{\d x}+2(n-1)BP=0.
\ee

\end{document}